\documentclass[12pt,reqno]{amsart}

\usepackage{epsfig}

\renewcommand{\Re}{{\operatorname{Re\,}}}
\newcommand{\inv}{^{-1}}
\newcommand{\szego}{Szeg\"o\ }
\newcommand{\sm}{\setminus}

\newcommand{\wt}{\widetilde}
\newcommand{\wh}{\widehat}
\newcommand{\PP}{{\mathbb P}}
\newcommand{\R}{{\mathbb R}}
\newcommand{\C}{{\mathbb C}}

\newcommand{\CP}{\C\PP}
\renewcommand{\d}{\partial}
\newcommand{\dbar}{\bar\partial}
\newcommand{\ddbar}{\partial\dbar}

\newcommand{\E}{{\mathbf E}}

\newcommand{\half}{{\frac{1}{2}}}
\newcommand{\vol}{{\operatorname{Vol}}}

\newcommand{\SU}{{\operatorname{SU}}}

\newcommand{\supp}{{\operatorname{Supp\,}}}
\renewcommand{\phi}{\varphi}
\newcommand{\eqd}{\buildrel {\operatorname{def}}\over =}

\newcommand{\ccal}{\mathcal{C}}
\newcommand{\dcal}{\mathcal{D}}

\newcommand{\hcal}{\mathcal{H}}

\newcommand{\lcal}{\mathcal{L}}
\newcommand{\pcal}{\mathcal{P}}

\newcommand{\scal}{\mathcal{S}}

\newcommand{\al}{\alpha}

\newcommand{\ga}{\gamma}

\newcommand{\La}{\Lambda}
\newcommand{\la}{\lambda}
\newcommand{\ep}{\varepsilon}
\newcommand{\de}{\delta}

\newcommand{\om}{\omega}

\newtheorem{theo}{{\sc Theorem}}

\newtheorem{lem}[theo]{{\sc Lemma}}
\newtheorem{prop}[theo]{{\sc Proposition}}

\newenvironment{rem}{\medskip\noindent{\it Remark:\/} }{\medskip}

\title
{Equilibrium distribution of zeros of random polynomials }

\author{Bernard Shiffman}
\author{Steve Zelditch}
\address{Department of Mathematics, Johns Hopkins University, Baltimore,
MD
21218, USA}
\email{bshiffman@jhu.edu, szelditch@jhu.edu}

\thanks{Research partially supported by NSF grants
\#DMS-0100474 (first author) and \#DMS-0071358 (second author).}

\date{June 17, 2002}

\begin{document}

\begin{abstract}
We consider ensembles  of random polynomials of the form $p(z) =
\sum_{j = 1}^N a_j P_j$ where $\{a_j\}$ are independent complex
normal random variables and where $\{P_j\}$ are the orthonormal
polynomials on the boundary  of  a bounded simply connected
analytic plane domain $\Omega \subset \C$ relative to an analytic
weight $\rho(z) |dz|$ . In the simplest case where $\Omega$ is the
unit disk  and $ \rho =  1$, so that   $P_j(z) = z^j$, it is known
that the average distribution of zeros is the uniform measure on
$S^1$.
 We  show that for any analytic $(\Omega, \rho)$,
    the   zeros of  random polynomials almost surely become
equidistributed relative to the equilibrium measure on $\d\Omega$
as $N\to\infty$. We further show that on the length scale of
$1/N$, the correlations have a universal scaling limit independent
of $(\Omega, \rho)$. 
\end{abstract}

\maketitle

\section*{Introduction}

A well-known result due to Hammersley \cite{H} (see also
Shepp-Vanderbei \cite{ SV}) states that  the zeros of random
complex `Kac'  polynomials
\begin{equation} \label{POLY} f(z) = \sum_{j = 0}^N a_j z^j,\;\;\; z \in \C
\end{equation} tend to concentrate on the unit circle $S^1 = \{|z| = 1\}$ as
$N \to
\infty$ when the coefficients $a_j$ are
 independent complex Gaussian random variables of mean $0$ and variance $1$:
\begin{equation} \label{SONE} {\bf E}(a_j) = 0,
\quad  {\bf E}(a_j \bar{a}_k) =
\delta_{j k}, \quad  {\bf E}(a_j {a}_k) =0. \end{equation} To be precise,
(\ref{SONE}) defines a Gaussian probability measure $\gamma^N_{S^1}$ on the
space ${\mathcal P}_N$ of polynomials of degree $\le N$ on $\C$,
studied long ago  by Littlewood-Offord, Erdos-Turan, and in
particular by Kac \cite{K1, K2}.  We denote expected values
relative to $\ga^N_{S^1}$ by ${\bf E}_{S^1}^N$. We also define
the {\it normalized distribution of zeros\/} of $f$ to be
the probability measure
\begin{equation}\label{nzd}\wt Z^N_{f}:=  \frac{1}{N} \sum_{f(z) = 0}
\delta_z\;.\end{equation} Then we may formulate the  concentration of zeros
on
$S^1$ as the weak limit formula:
\begin{equation} {\bf E}_{\gamma^N_{S^1}}(\wt Z^N_f) \to \delta_{S^1} \quad
\mbox{as }\
N
\to
\infty\;, \quad (\delta_{S^1},\phi) := \frac 1{2\pi}\int_{S^1} \phi(e^{i
\theta})\, d
\theta.
\label{hammer}\end{equation}

The first purpose of this note is to  generalize the expected
equidistribution result (\ref{hammer}) in a suitable sense to  any
closed analytic curve $\partial \Omega$ in $\C$ which bounds a
simply connected plane domain $\Omega$.  Our method is new even in
the case of $S^1$ and gives sharper results than
those of \cite{H}.  We further prove that the equidistribution
result is {\it self-averaging} in the sense that $\wt Z_{f_N}$ is almost
surely
asymptotic to the expected distribution. Our second purpose is to show that
the
correlations between zeros have a universal scaling limit independent of
$\Omega,
\rho$ on the length scale of $1/N$ and to determine the properties of the
limit
pair correlation function.

Our starting point is to reinterpret the Gaussian measure
(\ref{SONE}) on ${\mathcal P}_N$ as the one induced by the inner
product
\begin{equation} \label{circle}\left\langle f,\bar g\right \rangle_{S^1} =
\frac 1{2\pi}
\int_{S^1} f\bar g\,|dz| \end{equation}
on ${\mathcal P}_N$.
In general, any inner product
$$\left\langle f,\bar g\right \rangle_\mu = \int_{\C} f\bar g\,d
\mu$$ on
${\mathcal P}_N$ induces a Gaussian measure on ${\mathcal P}_N$ as follows:
Let
$\{P_j\}$ denote an orthonormal basis of ${\mathcal P}_N$ relative to
$\left\langle , \right \rangle_\mu$ and write any polynomial in this basis:
\begin{equation}\label{gaussian}f(z) = \sum_{j = 0}^N a_j
P_j. \end{equation} The Gaussian measure is defined by the
condition that the coefficients  $a_j$ in this basis are i.i.d.\
complex Gaussian variables with mean zero and variance one; i.e.
it equals $\pi^{-N}  e^{-|a|^2} da$ in terms of the
coefficients.

The choice of the  Gaussian measure $\gamma_{S^1}^N$ in
(\ref{SONE}) results in
 the zeros of the random polynomials becoming equidistributed relative to
 $\delta_{S^1}$. As is well-known, $\delta_{S^1}$ is the equilibrium measure
of $S^1$ (or of the closed unit disk). Recall that the equilibrium measure
of
a compact set $K$ is the unique probability measure $\nu_K$ which minimizes
the energy
$$E(\mu) = -\int_K \int_K \log |z - w| \,d\mu(z)\, d \mu(w). $$
We claim that the same equidistribution  result holds if we replace the unit
disk
$U=\{|z|<1\}$ by any simply connected, bounded plane domain $\Omega$ with
real-analytic boundary.  That is, we replace the inner product
(\ref{circle}) on $S^1$ with an inner
product on $\d\Omega$ of the form
\begin{equation}\label{bdry}\langle f, \bar g \rangle_{\partial \Omega,
\rho} : =
\int_{\partial
\Omega} f(z)
\overline{g(z)}
\rho(z)\, |dz|\;,\qquad \rho \in \ccal^\omega(\d\Omega)\;. \end{equation} We
denote by $\gamma^N_{\partial \Omega, \rho}$ the Gaussian measure
induced by this inner product on ${\mathcal P}_N$ as in
(\ref{gaussian}), and we denote the expectation relative to the
ensemble $(\pcal_N,\,\gamma^N_{\partial \Omega, \rho})$  by ${\bf
E}^N_{\partial \Omega, \rho}$. As in the $S^1$ case, we  obtain an
asymptotic formula
for the expected normalized distribution of zeros ${\bf E}^N_{\partial
\Omega,
\rho}(\wt Z^N_f)$:

\begin{theo}\label{OMEGA} Suppose that $\Omega$ is a simply-connected
bounded
$\ccal^\om$ domain and $\rho$ is a positive $\ccal^\om$ density on $\d
\Omega$.
Then with the above notation,
$${\bf E}^N_{\partial \Omega, \rho} (\wt Z^N_f)
 = \nu_\Omega
+O\left(1/N \right)\;,$$
where $\nu_\Omega$ is the equilibrium measure of $\bar \Omega$. \end{theo}

Here we
use the expression $O(1/N)$ to mean a distribution $u_N\in\dcal'(\C)$ such
that
$|(u_N,\phi)|\le C_\phi/N$  for
all $\phi\in\dcal(\C)$, where $C_\phi$ independent of $N$. Note that the
limit
distribution of zeros concentrates on $\partial \Omega=\supp(\nu_\Omega)$
and is
independent of the density $\rho$ used to define the inner product.

The equilibrium measure may be understood as the measure by which $N$
`electric charges'
in $\Omega$ distribute themselves in the limit $N \to \infty$.  Thus, at
least on average, the zeros of random
polynomials in the ensemble $(\pcal_N,\,\gamma^N_{\partial \Omega, \rho})$
behave like
electric charges. This suggests that the same equidistribution result should
hold if we
orthonormalize the polynomials within $\Omega$. We restrict to Lebesgue
measure in
$\Omega$ and define the inner product
\begin{equation}\label{int}\langle f, \bar g \rangle_\Omega : =
\int_{\Omega} f(z)
\overline{g(z)}
\,dx\,dy\qquad (z=x+iy)\;,\end{equation} which induces the Gaussian measure
$\ga_\Omega^N$
on
$\pcal_N$.

\begin{theo}\label{interior} Suppose that $\Omega$ is a simply-connected
bounded
$\ccal^\om$ domain and let $\E_\Omega^N$ denote the expectation with respect
to the
Gaussian measure $\ga_\Omega^N$ on
$\pcal_N$. Then
$${\bf E}^N_{\Omega} (\wt Z^N_f,\phi)
 = \nu_\Omega
+O\left(1/N\right)\;.$$
\end{theo}

As in \cite{SZ}, the expectation results of Theorems
\ref{OMEGA}--\ref{interior} can
be improved to give an almost sure limit distribution of zeros of random sequences
of polynomials of increasing degree:

\begin{theo}\label{as} Suppose that $\Omega$ is a simply-connected bounded
$\ccal^\om$ domain. Let $\ga_N$ be the Gaussian measure on $\pcal_N$ induced
by the
inner product (\ref{bdry}) or (\ref{int}), and let
$\mu:=\prod_{n=1}^\infty\ga_N$ be
the product probability measure on $\scal:=\prod_{n=1}^\infty\pcal_N$. Then
$\wt Z^N_{f_N}\to \nu_\Omega$ in the measure sense,  for
$\mu$-almost all sequences
$\{f_N\}\in\scal$; i.e.,
$$\frac 1N \sum_{\{z\in\C: f_N(z)=0\}} \phi(z) \to \int
_{\d\Omega}\phi\,d\nu_\Omega\qquad \forall \phi\in\ccal(\wh \C)\;,\qquad
\mu\mbox{-almost surely}\;.$$
\end{theo}

The principal  ingredients in the proofs of Theorems \ref{OMEGA} and
\ref{interior} are classical results of \szego \cite{Sz1} and of Carleman
\cite{Ca}
on the asymptotics of orthogonal polynomials normalized on the boundary
and on the interior of a domain,  respectively. We will recall these results
in \S
\ref{asympt}. The proofs of the two theorems are then essentially the same.

Theorems \ref{OMEGA}--\ref{as} say that in the limit the zeros tend to be
uniformly distributed along $\d\Omega$ (with respect to the equilibrium
measure).  However, as in our earlier work with P. Bleher \cite{BSZ1,BSZ},
the zeros do not behave as if they were thrown down at random; indeed, the
zeros are correlated.  To quantify this correlation, we consider as in
\cite{BSZ1,BSZ} the {\it $\ell$-point correlation functions\/}
\begin{equation}\label{corfn} K^{\ell N}_{\partial \Omega,
\rho}(z_1,\dots,z_\ell):=
\frac{{\bf E}^N_{\partial \Omega, \rho} (\overbrace{\wt Z^N_f\times \cdots\times
\wt Z^N_f}^{\ell})}{{\bf E}^N_{\partial \Omega, \rho} (\wt Z^N_f)\times
\cdots\times{\bf E}^N_{\partial
\Omega,
\rho} (\wt  Z^N_f)}\;,\qquad (z_1,\dots,z_\ell)\in \C_\ell\;,
\end{equation} where $$\C_\ell=\{(z_1,\dots,z_\ell)\in \C^\ell: z_j\neq z_k \
\mbox{for}\ j\ne k\}\;.$$ (We punch out the big diagonal in $\C^\ell$ since the
numerator in (\ref{corfn}) has a singular part there.) The correlation functions
$K^{\ell N}_{\Omega}(z,w)$ are similarly defined. 

In particular, the
{\it pair correlation
function\/}
$K^{2N}_{\partial
\Omega, \rho}(z,w)$ (or $K^{2N}_{\Omega}(z,w)$) can be interpreted as the
conditional probability of finding a zero at $w$ given that there is a zero
at
$z$, normalized by dividing by the unconditional probability of finding a zero at
$w$. For example, if $K^{2N}_{\partial \Omega, \rho}(z,w)=1$, then the
existence of a zero at
$z$ has no influence on the probability of finding a zero at $w$. This is
always the case in the limit as $N\to \infty$ if $z,w$ are fixed.  However,
we
have nontrivial correlations if the distance between points is
$O(1/N)$.  To describe this phenomenon, let $z_0\in \d\Omega$ be fixed, and
choose the complex coordinate $\zeta = \Phi(z)$, where $\Phi$ is the Riemann
mapping function mapping the exterior of $\Omega$ to the exterior of the
unit
disk, mapping $z_0$ to $1$, and  taking $\infty$ to itself. (This
exterior Riemann mapping function is a basic ingredient in the \szego and
Carleman asymptotics used in this paper.) We then have  universal scaling
limit zero correlation functions:

\begin{theo}\label{corthm}  There exist
universal functions
$K^{\ell\infty}:\C_\ell\to \R^+$ independent of $\Omega, z_0,\rho$ such that
$$\wh K^{\ell N}_{\partial \Omega, \rho}\left(1 +\frac {\zeta_1}
N, \dots, 1 +\frac {\zeta_\ell} N\right) \to 
K^{\ell\infty}(\zeta_1,\dots,\zeta_\ell)$$ as
$N\to\infty$, where $\wh K^{\ell N}_{\partial \Omega, \rho}=K^{\ell N}_{\partial
\Omega,
\rho}\circ \Phi\inv$ is
the correlation function written in terms of the complex coordinate $\zeta$
described above.  Similarly, $\wh K^{\ell N}_\Omega\left(1 +\frac {\zeta_1}
N, \dots, 1 +\frac {\zeta_\ell} N\right) \to 
K^{\ell\infty}(\zeta_1,\dots,\zeta_\ell)$.
\end{theo}
\noindent The universal scaling limit pair correlation function  $K^{2\infty}$ is given by
formula (\ref{formula}).

In our work with P. Bleher \cite{BSZ1,BSZ}, we showed that the zero correlation
functions can be given by  general (universal) formulas involving only the \szego
kernel and its first and second derivatives.  In  \S \ref{s-cor}, we study the
partial \szego kernel
\begin{equation}\label{pszego} S_N(z,w) =  \sum_{k = 0}^{N} P_k(z)
\overline{P_k(w)}\;,\end{equation} which gives the orthogonal projection
onto the span of the first $N$ orthogonal polynomials associated
to the inner product (\ref{bdry}). We show (Proposition~\ref{limszego}) that the
scaling asymptotics of $S_N$ have  the form
\begin{equation*} \frac 1N\wh S_N\left(1+\frac{\zeta_1}N, 1+\frac{\zeta_2}N\right)
 \to C_{\Omega,\rho,z_0}\;G(\zeta_1+\bar\zeta_2)\;,\end{equation*} 
where $\wh S_N = S_N\circ\Phi\inv$ and
\begin{equation}\label{G}G(z)=\frac{e^z-1}{z}\;.\end{equation} 

It is natural to separately consider the  `tangential' scaling
asymptotics along the boundary and the `normal' scaling
asymptotics orthogonal to the boundary. In the tangential case, we
set $\zeta_1 = i \theta, \zeta_2 = 0$ and obtain the tangential
scaling limit kernel $G(e^{i \theta}) = e^{i \theta/2} (\frac{\sin
\theta/2}{\theta/2}).$ This modified sine-kernel is reminiscent of
the scaling asymptotics of the projection kernels ${\mathcal K}_N$
for the  orthogonal polynomials which occur in the theory of
random matrices (see \cite{D}, Theorem 8.16). However, our actual
result is somewhat different and the methods have little in
common.

\begin{figure}[ht]\centering
\epsfig{file=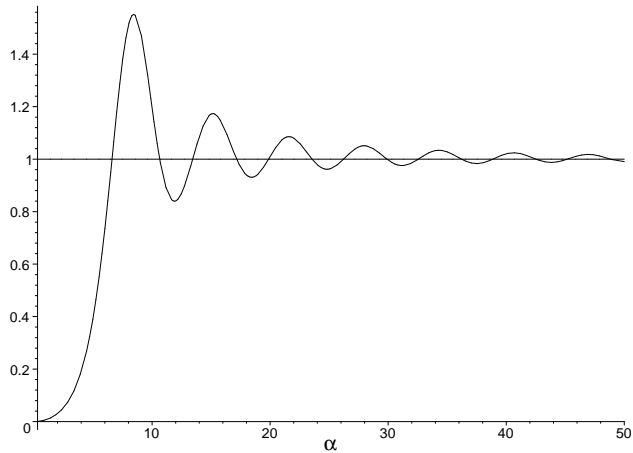,bb=72 122 540 544,height=2.5in}
\caption{The tangential  pair correlation function
$\kappa^{\rm T}(\al)$}\label{cor-graph}
\end{figure}

The tangential scaling limit for pair
correlations between zeros,
$$\kappa^{\rm T}(\al):=K^{2\infty}(0,{i\al})= \lim_{N\to\infty}\wh
K^{2N}_{\partial \Omega, \rho}\left(1,1+{\textstyle \frac{i\al}N}\right)=
\lim_{N\to\infty}\wh K^{2N}_{\partial \Omega, \rho}\left(1,e^{i\al/N}\right)\;,$$ 
measures the probability (density) of finding a pair of zeros in
small disks around two points on $\partial \Omega$ in terms of the ($\frac
1N$-scaled) angular distance $\al$ between them. As is illustrated by the  graph
of
$\kappa^{\rm T}$ (drawn using Maple$^{\rm TM}$) in Figure~\ref{cor-graph}, in
tangential directions the  zeros `repel' as $\al \to 0$  and are oscillatory with
$1/\al^2$ decay over long scaled distances, very much as for correlations
between eigenvalues of Gaussian random Hermitian matrices.

\begin{figure}[ht]\centering
\epsfig{file=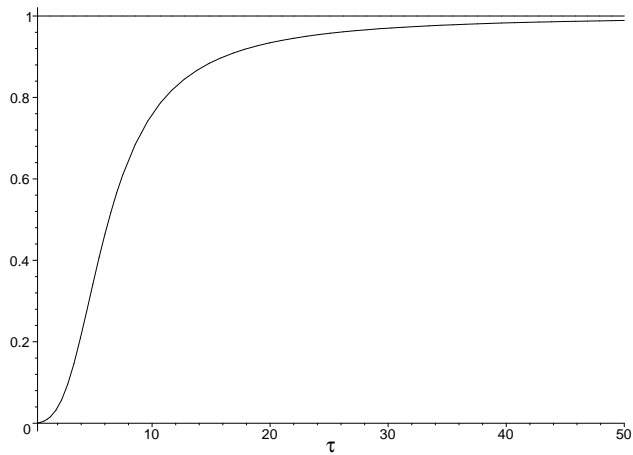,bb=72 122 540 544,height=2.5in}
\caption{The normal  pair correlation function
$\kappa^\bot(\tau)$}\label{cor-tran}
\end{figure}

In the normal directions, we obtain the scaling limit:
$$\kappa^{\bot} (\tau):=K^{2\infty}(0,\tau)= \lim_{N\to\infty}\wh
K^{2N}_{\partial \Omega, \rho}\left(1, 1+\frac\tau N\right)\;.$$ The graph
of $\kappa^\bot$ in Figure~\ref{cor-tran} is reminiscent of the
Bleher-Di correlation function \cite{BD} for real Kac polynomials.
Although the normal correlation is not oscillatory, we again have
zero repulsion and $1/\tau^2$ decay.

We  end the introduction with some remarks on and comparison to
prior results in this area. Real Kac polynomials, i.e. polynomials
as in (\ref{POLY}) but with real Gaussian coefficients of mean
zero and variance one, were shown by  Kac \cite{ K1, K2} to have
on average around  $\log N$ real zeros. The same was proved for
more general Kac polynomials in \cite{IM}. As mentioned above,
Hammersley \cite{H} obtained the concentration of zeros of complex Kac
polynomials on the unit circle in a somewhat less precise sense
than the one we state here. Later, Shepp and Vanderbei \cite{SV}
showed that the complex zeros of random real Kac polynomials also
tend to concentrate on the unit circle, explaining in a
qualitative way why so few zeros are real. Other classical studies
of random real polynomials are given in \cite{EO, LO, ET, EK}. The
probability of a
real Kac polynomial having no real zeros is estimated in \cite{DPSZ}.
Other results on hole probabilities for complex zeros are given in \cite{So}.

Recently, a number of  results were obtained on the
distribution of and correlations between zeros of Gaussian random
polynomials in one variable (and more general holomorphic
sections)  in which the inner product comes from a hermitian
metric $h$ on a positive line bundle $L \to M$ over a Riemann
surface $M$ (see \cite{Han, BR, SZ, BSZ1}. For simplicity, let
us mention only the case of $\SU(2)$ polynomials, in which $M =
\CP^1, L = {\mathcal O}(N), $ and where $h = h_{FS}^N$ is the
$N^{\rm th}$ power of the Fubini-Study metric. It is obvious from the
$\SU(2)$ symmetry of the expected distribution of zeros in this
case that it must equal the standard (Fubini-Study) area form, in
sharp contrast to the very singular equidistribution result we
find in Theorems \ref{OMEGA}--\ref{as}. Moreover, although the
pair correlation function for $\SU(2)$ polynomials given by \cite{Han} (see
\cite{BSZ1} for a generalization) exhibits the same repulsive behavior near 0, its
scaling limit occurs on the length scale of $1/ \sqrt{N}$ rather than
$1/N.$  Thus, the correlations are quite incomparable in the two
types of ensembles.

To explain the strong differences in the results, it should be
observed that  the ensemble $({\mathcal P}_N, \gamma_N)$ in the
$\SU(2)$ case comes from an inner product $\langle,\rangle_N$ on
${\mathcal P}_N$ which changes with each $N$, according to which
$\|z^k\|^2_N = {N\choose k}$. These inner products
$\langle,\rangle_N$ derive from a single inner product on CR
(Cauchy-Riemann) functions on the boundary of the unit ball in
$\C^2$ rather than on polynomials in one complex variable. This
suggests that  polynomials in the $\SU(2)$ ensemble should be
thought of as homogeneous polynomials $\sum_{j = 1}^N a_j z_0^{N-j}
z_1^j$ in two variables and that their equidistribution law and
correlations reflect the extra dimensionality. Because of the
concentration of zeros on curves, the
 ensembles studied here share a  one -dimensionality with the ensembles of
GUE and CUE random matrices
and with the observed correlations between  zeros of the Riemann
zeta function.

We note that in random matrix theory, orthogonal polynomials on
$\R$ with respect to weights $\rho_N = e^{- N V(x)} dx$ depending
on the degree $N$ are quite important (see \cite{D}). This suggests that
ensembles $({\mathcal P}_N, \gamma_N)$ of random polynomials on
curves $\partial \Omega$  where the inner product is given by a
 weight of the form $\rho_N = \rho(z)^N\,|dz|$ might be interesting. Such ensembles
have features in common with both the
 $\SU(2)$ ensemble (and its generalization to other Riemann surfaces) and
with the
ensembles of this paper. How are
 the zeros distributed and correlated? Do the zeros concentrate on
 $\partial \Omega$ or are they diffuse?
Our methods do not seem to adapt in a simple way
 to such ensembles, and perhaps they belong to different
 universality classes than the ones studied here.

Another interesting direction would be to  generalize our results
to higher dimensions.  Clearly the methods we use, based on the
Riemann mapping theorem,  have no simple generalization.

\section{Asymptotics of orthogonal polynomials}\label{asympt}

The proofs of Theorems~\ref{OMEGA}-\ref{interior} are   based on asymptotic
properties of  orthogonal  polynomials associated  to $\Omega$, and  to
their relations with the interior and
exterior Riemann mapping functions. These imply asymptotic properties of
partial
\szego and Bergman kernels for $\Omega.$
We recall the results we need in this section.
Classical references are \cite{Sz1, Ca} while more recent references (with
more
precise results) are in \cite{SL, Su}.

\subsection{\szego kernel and orthogonal polynomials}

Let $\Omega \subset \C$ be a smooth bounded domain. The \szego
kernel of $\Omega$ with respect to a measure $\rho |dz|$ on $\d\Omega$  is the
orthogonal
projection
\begin{equation} S: \lcal^2(\partial \Omega, \rho |dz|) \to \hcal^2(\partial
\Omega, \rho |dz| ) \end{equation} onto the Hardy space of  boundary values of
holomorphic functions in $\Omega$ which belong to $\lcal^2(\partial
\Omega, |dz|)$. The Schwartz kernel of $S$ is denoted $S(z,w)$.
According to \cite{B}, Theorem 24.3, $S(z,w)$ admits an analytic
continuation
(holomorphic in $z$ and anti-holomorphic in $w$)
to $\overline{\Omega} \times \overline{\Omega} \sm \Delta$, where $\Delta $
is
the diagonal in $\partial \Omega \times \partial \Omega.$ We will refer to
this
as the regularity theorem.

 Let
$$\{P_{ j} (z) = a_{ j0} + a_{j 1} z + \dots +a_{ j j} z^j\}$$
be the orthonormal basis of orthogonal polynomials for $\lcal^2(\partial
\Omega,
\rho |dz|)$ obtained by applying Gram-Schmidt to
$\{1,z,z^2,\dots,z^j,\dots\}$. Since $S$
is an orthogonal projection,we may express it in terms of any orthonormal
basis. Hence,
we have
\begin{equation} S(z,w) = \sum_{k = 0}^{\infty} P_k(z)
\overline{P_k(w)},\;\; (z,w) \in \overline{\Omega} \times \overline{\Omega}
\end{equation}
By the regularity theorem, one has that $S(z,z) < \infty$ for $z\in\Omega $,
and
thus
$P_N \to 0$ on $\Omega$. Hence,
$$S_N(z,z) \to S(z,z),\;\; \mbox{uniformly on compact subsets of} \; \Omega
,$$
where
$S_N(z,w) =  \sum_{k = 0}^{N} P_k(z)
\overline{P_k(w)}$ is the partial \szego
kernel given in (\ref{pszego}).

In  the  case of the unit disk $U$ (with the weight $\rho\equiv 1$), one
has:
$$S^U(z,w) = \frac{1}{2 \pi (1 - z \bar{w})}.$$
We observe that $S(z,w)$ admits a meromorphic continuation to the exterior
as well as the interior of $\Omega.$ However, the behavior is quite
different
with regard to the partial \szego kernels. Indeed,
The orthogonal polynomials are $P_k(z) = z^k$, hence
\begin{equation}\label{SNU}|P_N(z)|^2 = |z|^{2N},\quad S_N^U(z, z) = \sum_{k
= 0}^{N- 1}
|z|^{2k} =
\frac{1 - |z|^{2N}}{1 - |z|^2}\;.\end{equation} Clearly, $S_N(z,z) \to
\infty$ at an
exponential rate in the exterior of $\Omega.$

\subsection{Partial \szego kernels and the exterior Riemann mapping
function}

We will need the behavior of $S_N(z,z)$ in the exterior of
$\Omega$. The relevant result was proved by \szego \cite{Sz1,Sz2}.

Let $\Omega \subset \C$ denote a simply connected bounded plane domain with
$C^{\infty}$ boundary $\partial \Omega$. The exterior domain $\wh\C \sm
\overline\Omega$ is also a simply connected domain in $\wh{\C}$ and we
denote by
\begin{equation} \Phi: \wh\C \sm \Omega \to \wh \C\sm U\;,\qquad
\Phi(z) = cz + c_0 + c_1 z^{-1} + \cdots \end{equation} the
(unique) exterior  Riemann mapping function with $\Phi(\infty) =
\infty$, $\Phi'(\infty)\in\R^+$. We recall that the equilibrium
measure $\nu_{\Omega}$ of $\Omega$ is given by
\begin{equation}\label{equil}\nu_{\Omega} = \Phi^*\de_{S^1}\;;
\qquad \mbox{i.e., }\ \int\phi\,d\nu_\Omega= \frac {1}{2\pi}
\int_0^{2\pi}\phi\circ \Phi\inv(e^{i\theta})\,d\theta\;.\end{equation}

\begin{rem} If we fix
a point $x_0 \in \Omega$, we may define as well the interior Riemann mapping
function $f:\Omega\to U$ with $f(x_0) = 0,\ f'(x_0)\in\R^+$.  The functions
$f|_{\d\Omega}$ and
$\Phi|_{\d\Omega}$ are generally different.\end{rem}

Let us assume that $\partial \Omega \in \ccal^\omega $. One can define a
(unique)
nonvanishing holomorphic function $\Delta_e$ on $\C \sm \Omega$ satisfying:
$$|\Delta_e|^2_{\d\Omega}=\rho\;.$$ (See \cite[\S\S 10.2,16.1]{Sz2}.) In the
following, $L=\mbox{Length}(\d\Omega)$, and
$$T_{\epsilon}(\partial \Omega) = \{x \in \C: d(x, \partial \Omega) <
\epsilon\}.$$

\begin{theo} \label{SZEGO} \cite[Theorems 16.4--16.5]{Sz2} Let $\{P_n\}$
denote
the orthonormal polynomials relative to $(\partial \Omega, \rho |dz|)$, where
$\Omega,\,
\rho$ are as in Theorem~\ref{OMEGA}. Normalize
$P_n$ so that
its leading coefficient is positive.  Then there exist $\ep>0$ and $0<\de<1$
such
that for $z\in (\C \sm
\Omega)
\cup T_{\ep}(\partial
\Omega)$, we have
$$P_n(z) = \left(\frac{L}{2\pi}\right)^{1/2} \Delta_e(z)^{-1} \Phi'(z)^{1/2}
\Phi(z)^n + O(\de^n)\;. $$
\end{theo}

The proof is based on the properties of Faber polynomials $F_n$  associated
to
$G_n(z):=\left(\frac{L}{2\pi}\right)^{1/2}\Delta_e(z)^{-1} \Phi'(z)^{1/2}
\Phi(z)^n $; i.e., $F_n(z)$ is the polynomial part of the Laurent expansion
of
$G_n(z)$ about $z=\infty$. Suppose that the Riemann mapping function
$\Phi\inv$ extends analytically to $\{|z| > r\}$ where $r<1$. The main
estimates are:
$$P_n = F_n + O(n^{1/2} r^{n/2}) = G_n+ O(n^{1/2} r^{n/2})\;. $$ This
shows that the Faber polynomials (with weights) are asymptotic to the
orthogonal
polynomials (with weights).

\subsection{Bergman kernel and orthogonal polynomials}
In the proof of  Theorem \ref{interior} on the expected
distribution of zeros, where the inner product is chosen to be
$\langle f,\bar g\rangle_\Omega = \int_{\Omega} f\bar g \,dx\,dy$,
the  role of the \szego kernel is  played by the Bergman kernel,
i.e.  the orthogonal
 projection from $\lcal^2(\Omega)$ onto the subspace $\hcal^2(\Omega)$
spanned
by the
$\lcal^2$
 holomorphic functions.

As before, we assume  $\Omega \subset \C$ is a simply
connected, bounded plane domain with $\ccal^\omega$ boundary.
We denote by
$$\{P_{ j} (z) = a_{ j0} + a_{j 1} z + \dots +a_{ j j} z^j\}$$
 the orthonormal basis of orthogonal polynomials for $\lcal^2(\Omega, dxdy)$
with
positive leading coefficient. The Bergman kernel
 may be expressed in terms of the
 orthogonal polynomials by:
\begin{equation} B(z,w) = \sum_{k = 0}^{\infty} P_k(z)
\overline{P_k(w)}\;,\quad (z,w) \in {\Omega} \times {\Omega}.
\end{equation}
We let $B_N(z,w) = \sum_{k = 0}^N P_k(z)
\overline{P_k(w)}$; as in the case of the \szego kernel, we have
$$B_n(z,w)\to B(z,w)\;,\qquad B(z,z)>0\;.$$

 Orthogonal  polynomials over a domain  were studied by Carleman after
Szeg\"o's work on polynomials orthogonal on $\partial \Omega$. A key element
in our
proof is the following asymptotic formula for $P_N(z)$:

\begin{theo} {\rm (Carleman \cite{Ca}; cf.\ \cite{SL}, Theorems 1--2, p.
290)}  Let
$\{P_n\}$ denote the orthonormal polynomials relative to $(\Omega, dx\,dy)$,
where
$\Omega$ is a simply-connected bounded $\ccal^\omega$ domain.  Normalize
$P_n$ so that
its leading coefficient is positive.
Then there exist $\ep>0$ and $0<\de<1$ such that:
$$P_n(z) = \left\{\begin{array}{ll} \left(\frac{N + 1}{\pi}\right)^{1/2}
\Phi'(z)
\Phi(z)^n + O(\de^n),\qquad & z\in T_{\ep}(\partial
\Omega)\cap\Omega\\[10pt]
 \left(\frac{N + 1}{\pi}\right)^{1/2}
\Phi'(z)
\Phi(z)^n[1 + O(\de^n)],\qquad & z\in \C \sm \Omega\end{array}\right.$$
\label{Carle}\end{theo}

\section{Proof of the asymptotic formulas}

We now prove Theorems \ref{OMEGA}--\ref{interior}. The proofs are basically
the same in the boundary and interior cases, granted the theorems of \szego
and Carleman. We therefore only give details in the boundary case.

We first relate zero distributions to \szego kernels, as in our previous
papers
\cite{SZ,SZ2}.

\subsection{Zero distributions} For convenience, we let
$$Z^N_{f}=  \sum_{f(z) = 0}
\delta_z$$ denote the zero distribution of a polynomial $f$ so that
recalling
(\ref{nzd}), we have $\wt Z_f=\frac 1N Z_f$. (Strictly speaking, we should
count
the zeros with multiplicities, but polynomials in $\pcal_N$ almost surely
have only
simple zeros.) The expected zero distribution has a simple expression:

\begin{prop} \label{poincare} We have
$$ \E^N_{\partial \Omega, \rho} (Z_f) = \frac{\sqrt{-1}}{2\pi}\ddbar \log
S_N(z, z). $$ \end{prop}

This proposition is a special case of Proposition~4.1 in \cite{SZ2} (which
is a
variation of a result in our earlier work
\cite{SZ}).  For completeness, we give the proof here:   We first note that
since
$$ Z_f = \frac{\sqrt{-1}}{2\pi}\ddbar \log |f|^2, $$
we have
$${\bf E}^N_{\partial \Omega, \rho} ( Z_f) =
\frac{\sqrt{-1}}{2\pi}\ddbar{\bf E}^N_{\partial \Omega, \rho}\left( \log
|f|^2\right)\;.$$ To calculate the expectation, we write $f$ in
terms of the orthonormal basis $\{P_j\}$ of ${\mathcal P}_N$:
$$f(z) = \sum_{j = 0}^N a_j P_j(z)= \langle a, p(z) \rangle\;, $$ where
$a=(a_0,\dots,
a_N), $ $P=(P_0,\dots, P_N)$. Then,
$${\bf E}^N_{\partial \Omega, \rho} ( Z_f) = \frac{\sqrt{-1}}{\pi}\ddbar
\int_{\C^{N+1}}
\log |\langle a, P(z) \rangle|\, \frac 1{\pi^{N+1}} e^{-\|a\|^2}\,
da.$$ We write
$$ P(z) = \|P(z)\| u(z), \;\ \|P(z)\|^2 = \sum_{j = 0}^N |P_j(z)|^2 =
S_N(z,z),\;\ \|u(z)\|=1\;.$$ Then,
$$\log |\langle a, P(z) \rangle| = \log \|P(z)\| + \log |\langle a, u(z)
\rangle|\;.$$ We observe that
$$\int_{\C^{N+1}}
\log |\langle a, u(z) \rangle |\, e^{-\|a\|^2}\, da = \mbox{constant}$$
since for each $z$ we may apply a unitary coordinate change so that $u(z)
=(1,0,\dots,0)$. Hence the derivative equals zero, and we have
\begin{equation*} {\bf E}^N_{\partial \Omega, \rho} ( Z_f) =
\frac{\sqrt{-1}}{\pi}\ddbar
\log \|P(z)\|= \frac{\sqrt{-1}}{2\pi}\ddbar \log S_N(z,z)\;.
\end{equation*} \qed

By exactly the same argument, we also have:

\begin{prop} \label{poincare2} We have
$$ \E^N_{\Omega} (Z_f) = \frac{\sqrt{-1}}{2\pi}\ddbar \log B_N(z, z). $$
\end{prop}

\medskip
\subsection{The circular ensemble}

We start with the fundamental case of the circle; in the next section we will
reduce the other cases to the circular case.  Before proving our precise result
for the circle (Proposition~\ref{circular}), we
first give a simple ad hoc argument that the expected
limit distribution of zeros in this case is the measure $\de_{S^1}$
given in ({\ref{hammer}): By (\ref{SNU}),
$$\frac{1}{N} \ddbar \log S_N(z,z) \sim \frac{1}{N} \ddbar \log  (1 -
|z|^{2N}).$$ Clearly, in any annulus $|z|\leq  r < 1$, $(1 - |z|^{2N}) \to
1$ rapidly
with its derivatives, and the limit equals zero. In any annulus $|z| \geq r
> 1$ we
may write  $(1 - |z|^{2N}) = |z|^{2N}(|z|^{-2N} - 1)$ and separate the
factors
after taking $\log$. The second again tends to zero rapidly, while the first
factor,
 $\log  |z|^{2N}$, is killed by $\ddbar$ (note that $z \not= 0$
in this part).  It follows that the limit
measure must be supported on $S^1$. Since it is SO$(2)$-invariant (radial),
and
since it is a probability measure, it
must be $\frac{1}{2\pi}d \theta$, which we henceforth denote by $\nu$.

We have the following explicit formula and asymptotics for the
circular
case:

\begin{prop} \label{circular}   Let $\nu=\frac{d\theta}{2\pi}$ denote Haar
measure on
$S^1$.  Then
$$\E_\nu^N(Z_f)= \left[\frac 1{(|z|^2-1)^2} -\frac{(N+1)^2
|z|^{2N}}{(|z|^{2N+2}-1)^2}\right] \frac{\sqrt{-1}}{2\pi}\,dz\wedge d\bar
z\;,$$
Furthermore, $\E_\nu^N(Z_f)= N\nu +O(1)$; i.e., for all test forms
$\phi\in\dcal(\C)$, we have
$$\E_\nu^N\left(\sum_{\{z:f(z)=0\}}\phi(z)\right)= \frac{N}
{2\pi}\int_0^{2\pi}\phi(e^{i\theta})\,d\theta +O(1)\,.$$  In
particular, $ \E_\nu^N(\wt Z^N_f) \to\nu$ in $\dcal'(\C)$.
\end{prop}

The formula for $\E_\nu^N( Z_f)$ in Theorem \ref{circular} agrees with the
one given
by Hammersley \cite{H} (see also \cite{EK,SV}).

\begin{proof} By Proposition \ref{poincare},
\begin{equation} \E_\nu^N(Z_f)= \frac{\sqrt{-1}}{2\pi}\ddbar \log
\frac{1-|z|^{2N+2}}{1-|z|^{2}}\;.\label{U}\end{equation}

We write $\zeta={\rho +i\theta}=2\log z$, and we let
\begin{equation}\label{hN}h_N(\rho)=S_{N}^U(e^{\half\rho},e^{\half\rho})=
\sum_{n=0}^N e^{n\rho}=
\frac{1-e^{(N+1)\rho}}{1-e^\rho}\;.\end{equation}
Then \begin{equation}\label {Enu}\E_\nu^N( Z_f) = \frac{i}{8\pi}(\log
h_N)''(\rho)
\,d\zeta\,d\bar
\zeta\;.\end{equation} We let \begin{equation}\label{comp} g_{N}:= (\log
h_N)''
=\frac{e^{-\rho}}{(1-e^{-\rho})^2}
-(N+1)^2\frac{e^{-(N+1)\rho}}{(1-e^{-(N+1)\rho})^2} \;.\end{equation}
The desired formula for $\E_\nu^N( Z_f)$ follows by substituting
$d\zeta=\frac 2z
dz$ into (\ref{Enu})--(\ref{comp}).

We easily see that $g_N>0$ (which also follows from the strict
subharmonicity
of\break
 $\log S_N^U(z,z)= \log (\sum_{n=0}^N |z|^{2n})$) and that
 $g_N(-\rho)=g_N(\rho)$.
Now suppose that $\psi\in\ccal^\infty(\R)\cap \lcal^\infty(\R)$ such that
$\psi(0)=0$.  We claim that
\begin{equation}\label{gN}\int g_N \psi =O(1) \;.\end{equation}
To verify (\ref{gN}), we break up the integral into 3 pieces: $\rho\le
-1,\  -1\le \rho\le 1$, and $ \rho\ge 1$.
First, we have
$$\int_{-\infty}^{-1} g_N = \int^{+\infty}_{1} g_N=
\frac 1 {e-1}-\frac {N+1}{e^{N+1}-1} \le \frac 1 {e-1}
\;,$$ and thus $$\left|\int_{-\infty}^{-1}
g_N\psi\right|+
\left|\int^{+\infty}_{1} g_N\psi\right|=O(1)\;.$$

To estimate the integral over $[-1,1]$, we write $\psi(\rho)= a\rho
+O(\rho^2)$ and
we  note that
$$g_N(\rho)\le
\frac{1}{(1-e^{-\rho})^2}\le \frac 4 {\rho^2} \qquad\mbox{for }\ 0\le\rho\le
1\;.$$
Since $g_N$ is even, we then have
$$\left|\int^1_{-1} g_N\psi\right| \le
C\int^1_0 \rho^2g_N(\rho)\,d\rho
\le 4C\;,$$
which completes  the  proof of (\ref{gN}).

Now let $\phi\in\dcal(\R)$.  An easy computation gives  $\int_\R g_N d\rho =
N$.
 (This also follows from the fact that the expected number of zeros of a
degree-$N$ polynomial is
$N$.)   Thus by (\ref{gN}),
$$\int_\R g_N \phi = \phi(0)\int_\R g_N + \int_\R g_N[\phi -\phi(0)] =
N\phi(0)
+O(1)\;.$$
\end{proof}

\begin{rem}
Recalling (\ref{hN}), we have
\begin{equation*}g_{N-1}(0)
=\frac{h_{N-1}(0)h''_{N-1}(0)-h'_{N-1}(0)^2}{h_{N-1}(0)^2}
=\frac{N\sum_{k=0}^{N-1}k^2-(\sum_{k=0}^{N-1}k)^2}{N^2} =
\frac{N^2-1}{12}\;.\end{equation*}
One can easily check that $g_N'(\rho)<0$ for $\rho>0$ and since $g_N$ is
even, we
have
$g_N(0)= \sup g_N$.
In fact, a computation using Maple$^{\rm TM}$ gives:
$$g_{N-1}=
\frac{{N}^{2}-1}{12} -\frac {N^4-1}{240}
 {\rho}^{2}+ {\frac {N^6-1}{6048
}}  {\rho}^{4} -{\frac {{N}^{8}-1}
{172800}}  {\rho}^{6}+ {
\frac {{N}^{10}-1}{5322240}}  {\rho}^{8} -{\frac {691 ({N}^{12}-1)}{
118879488000}}\rho^{10}\dots
$$
\end{rem}

\subsection{The general case}

We now prove Theorem \ref{OMEGA}.
We must show that \begin{equation} \label{omega-eq}
\frac{\sqrt{-1}}{2\pi} \ddbar \log S_N(z,z) = N\nu_\Omega
+O(1)\;.\end{equation}

We shall break up the proof of (\ref{omega-eq}) into two regions:  $\Omega$
and
a neighborhood $W_\ep$ of $\wh\C\sm\Omega$.  The estimate (\ref{omega-eq})
is obvious
on
$\Omega$. Indeed, by Theorem~16.3 in \cite{Sz2},
$S_N(z,z) \to S(z,z)$ uniformly in compact subsets of $\Omega $.
Furthermore,
$S(z,z)>0$ on $\Omega $ (see e.g., \cite[Th.~12.3]{B}), and hence $\log
S_N(z,z)
\to \log S(z,z)$.
Thus
$$\ddbar \log  S_N(z,z) = \ddbar\log  S(z,z) +o(1)=O(1)\quad \mbox{in }\
\dcal'(\Omega)\;.$$ (In fact,
$S_N(z,w)\to S(z,w)$ on $\Omega  \times \Omega $, and hence by the Cauchy
Integral
formula, we have normal convergence of all derivatives; it then follows that
$\ddbar \log  S_N(z,z) \to \ddbar \log  S(z,z)$
uniformly on compact subsets of $\Omega $.)

Next we verify (\ref{omega-eq}) on the domain
$$W_\ep:=(\wh\C\sm\Omega)\cup T_\ep(\d\Omega)\;,$$ where $\ep$ is chosen
sufficiently
small so that Theorem
\ref{SZEGO} holds and
$|\Phi| \ge \de'>\de$ on $W_\ep$. 

Let \begin{equation}\label{SN}A_N(z)=S_N(z,z)/\Phi_N^*S_N^U(z,z)\;\;.\end{equation}
We claim that there is a positive constant $C$ such that
\begin{equation}\label{AN}0<1/C\le A_N(z)\le C <+\infty\qquad \mbox{for }\ z\in
W_\ep\;.\end{equation} To verify (\ref{AN}), we let
\begin{equation}\label{psi}\psi(z)=
\left(\frac{L}{2\pi}\right)^{1/2}
\Delta_e(z)^{-1}
\Phi'(z)^{1/2}\;,\end{equation}   and we note that $\psi$ and $\frac{1}{\psi}$ are
bounded on
$W_\ep$.
Recalling that
$$\Phi_N^*S_N^U(z,z)= \sum_{n=0}^N |\Phi(z)|^{2n}\;,$$ it then follows from
Szeg\"o's
Theorem \ref{SZEGO} that $\sup_{z\in W_\ep}A_N(z) \le C
<+\infty$.

The lower bound for $A_N(z)$ follows from
\cite[\S
16.5]{Sz2}.  Alternately, let
$G_n(z):=\psi(z)
\Phi(z)^n $ as above.  Then by Theorem  \ref{SZEGO}, there exists $n_0\ge 2$
such
that $|P_n(z)|^2\ge \half |G_n(z)|^2$ for $z\in W_\ep,\ n\ge n_0$.  If
$\de'<\|\Phi(z)|\le 1$ we have
$|p_0(z)|^2=c\ge c'(|G_0(z)|^2+\cdots+ |G_{n_0}(z)|^2)$ and hence
\begin{equation}\label{SN2}S_N(z,z)\ge \min\{1/2,c'\}\sum_{n=0}^N |G_n(z)|^2
\ge c''
\Phi_N^*S_N^U(z,z)\;.\end{equation}
On the other hand, if $|\Phi(z)|\ge 1$ then the required
estimate follows  from
\begin{eqnarray*}S_N(z,z) &\ge &|P_{n_0}(z)|^2 +\cdots +|P_N(z)|^2\ \ge\
\half\left
[|G_{n_0}(z)|^2+\cdots +|G_N(z)|^2\right]\\ &\ge &
\frac 14 \left[|G_{0}(z)|^2+\cdots +|G_N(z)|^2\right]\ =\ \frac
14|\psi(z)|^2\Phi_N^*S_N^U(z,z)\;,\end{eqnarray*}
for $N\ge 2n_0$.

It follows from (\ref{AN}) that $$(\ddbar \log A_N,\phi) = 
(\log A_N,\ddbar \phi) = O(1)\;,\qquad \mbox{for }\ 
\phi\in
\dcal(W_\ep)\;.$$ 
Recalling (\ref{equil}), we then have \begin{eqnarray*} \frac{i}{2\pi} \ddbar \log
S_N(z,z) &=&
\Phi^*\left(\frac{i}{2\pi}\ddbar
\log S_N^U(z,z)\right) - \frac{i}{2\pi} \ddbar \log A_N(z)\\ &=& \Phi^*(N\nu +O(1))
+O(1) 
\ =\ N\nu_\Omega+O(1)\;.\end{eqnarray*}
\qed

\subsection{The interior case}

The proof of Theorem~\ref{interior} is exactly the same as in the boundary
case, using
Carleman's theorem in place of Szeg\"o's, and the Bergman kernel and
Proposition~\ref{poincare2} in place of the \szego kernel and
Proposition~\ref{poincare}. We omit the details.

We also note that the conclusion of Theorem \ref{interior} holds for inner
products on
$\Omega$ with certain analytic weights; see \cite{SL}.

\subsection{Proof of Theorem \ref{as}}  The proof of almost sure convergence
to the
average is exactly the same as the proof of the similar statement for
sections of
positive line bundles in \cite{SZ}.  We summarize the proof here:  Let
$$\omega_N=\left\{\begin{array}{llll}{\bf E}^N_{\partial \Omega,
\rho}(\wt Z^N_f) &=& \frac{\sqrt{-1}}{2\pi N}\ddbar \log S_N(z, z)\quad &
\mbox{for
the boundary case,}\\[10pt] {\bf E}^N_{\Omega}(\wt Z^N_f) &=&
\frac{\sqrt{-1}}{2\pi N}\ddbar
\log B_N(z, z)\quad & \mbox{for
the interior case.}\end{array}\right.$$
We then have the following variance estimate:
\begin{lem} \label{variance} Let
$\phi$ be any smooth test form.  Then $$\E\left(\left(\wt
Z^N_f-\omega_N,\phi\right)^2\right) = O(N^{-2}).$$ \end{lem}
Lemma \ref{variance} is  given as Lemma 3.3 in \cite{SZ} and has exactly the
same proof.  Continuing as in \cite{SZ}, by Theorem~\ref{OMEGA} or Theorem
\ref{interior}, it suffices to show that  $$(\wt
Z^N_{f_N}-\omega_N,\phi)\to 0 \quad\quad \mbox{\rm almost surely}\;.$$
Consider the random variables $X_N$ on $\scal$ given by: \begin{equation*}
X_N(\{f_N\}) =(\wt Z^N_{f_N}-\omega_N,\phi)^2\geq 0\;.\end{equation*} By
Lemma~\ref{variance},
$\int_{\scal} X_N\, d\mu = O(N^{-2})$, and therefore
$$\int_{\scal}\sum_{N=1}^\infty X_N\, d\mu = \sum_{N=1}^\infty \int_{\scal}
X_N\, d\mu <+\infty\;.$$ Hence, $X_N\to 0$ almost surely.  The conclusion
follows by
considering a countable $\ccal^0$-dense family of test forms.
\qed

\section{Pair correlations of zeros} \label{s-cor}
In our previous work with P. Bleher \cite{BSZ1,BSZ}, we derived the scaling limit
pair correlation functions for zeros of
$\SU(m+1)$-polynomials and showed that these correlations are `universal.'
In this section, we apply the general formulas from  \cite{BSZ1,BSZ} to
prove
Theorem~\ref{corthm} and to describe our universal pair correlation function
$K^{2\infty}$.

\subsection{The scaled \szego kernel} In 
\cite[\S 2.3]{BSZ1}, we gave a general formula for the $\ell$-point correlation of
zeros in terms of the Bergman-\szego kernel and its first and second derivatives. In
order to apply this formula to find the scaling limit correlations, we need the
following universal scaling limit \szego and Bergman kernels for plane domains:

\begin{prop} Let $\Omega$ be a simply-connected bounded $\ccal^\om$ domain and
let $\rho$ be a $\ccal^\om$ density.  
\begin{enumerate}
\item[\rm i)] Let $S_N$ be the orthogonal projection for
the inner product $\langle ,\rangle _{\d\Omega,\rho}$.  Then 
$$\frac 1N \wh S_N\left(1+\frac {\zeta_1}{N},1+\frac {\zeta_2}{N}\right) \to
|\psi(z_0)|^2 G(\zeta_1+\bar\zeta_2)\;,\qquad  G(z)=\frac{e^z-1}{z}\;,$$
where $\wh S_N= S_N\circ \Phi\inv$ is
the projection kernel written in terms of the complex coordinate
$\zeta=\Phi(z)$, $\psi$ is given by (\ref{psi}) and
$z_0=\Phi\inv(1)$.
\item[\rm ii)]  Let $B_N$ be the orthogonal projection for
the inner product $\langle ,\rangle _\Omega$.  Then 
$$\frac 1{N^2} \wh B_N\left(1+\frac {\zeta_1}{N},1+\frac {\zeta_2}{N}\right) \to
\frac 1\pi|\Phi'(z_0)|^2 G(\zeta_1+\bar\zeta_2)\;,$$ where $\wh B_N= B_N\circ
\Phi\inv$ is the projection kernel written in terms of the complex coordinate
$\zeta=\Phi(z)$ and
$z_0=\Phi\inv(1)$.
\end{enumerate}
The above limits are uniform when $|\zeta_1|+|\zeta_2|$ is bounded.
\label{limszego}\end{prop}

\begin{proof}  We begin as before with the case of the disk $U$ with density
$\rho=1$. Then
\begin{eqnarray} \frac 1N S_N^U(1+\zeta_1/N, 1+\zeta_2/N)
&=&\frac 1N \sum_{k=0}^N
\left(1+\zeta_1/N\right)^k\left(1+\bar\zeta_2/N\right)^k\nonumber\\ &=&
\frac{[(1+\zeta_1/N)(1+\bar\zeta_2/N)]^{\frac{N+1}{N}} -1}{\zeta_1+\bar\zeta_2
+\zeta_1\bar\zeta_2/N}\nonumber\\&\to &
\frac{e^{\zeta_1}e^{\bar\zeta_2}-1}{\zeta_1+\bar\zeta_2}  \ =\
G(\zeta_1+\bar\zeta_2)\;.\label{thedisk}\end{eqnarray}

We now consider the inner product $\langle ,\rangle _{\d\Omega,\rho}$.  By
Szeg\"o's Theorem \ref{SZEGO}, we can choose $\ep,\la>0$ such that
$$P_k(z)=\psi(z)[\Phi(z)^k+ q_k(z)]\,,\quad |q_k(z)|\le C_1e^{-\la k} \ \
\mbox{for }\ 1-\ep\le |\Phi(z)|\le 2\;,$$ and hence
\begin{eqnarray*}\frac {\wh S_N(1+{\textstyle \frac {\zeta_1}N}, 1+{\textstyle
\frac {\zeta_2}N})} {N \,\wh\psi(1+{\textstyle \frac
{\zeta_1}N})\overline{\wh\psi(1+{\textstyle \frac {\zeta_2}N})}} &=&
\frac 1N S_N^U(1+{\textstyle \frac {\zeta_1}N}, 1+{\textstyle \frac {\zeta_2}N}) 
\ +\ \frac 1N\sum_{k=0}^N(1+{\textstyle \frac
{\zeta_1}N})^k\,\overline{\hat q_k(1+{\textstyle \frac {\zeta_2}N})}\\ && +\ \frac
1N\sum_{k=0}^N\overline{(1+{\textstyle \frac {\zeta_2}N})^k}\,
{\hat q_k(1+{\textstyle
\frac {\zeta_1}N})} \ +\ \frac 1N\sum_{k=0}^N\hat q_k(1+{\textstyle \frac
{\zeta_1}N})
\overline{\hat q_k(1+{\textstyle \frac {\zeta_2}N})}\,.
\end{eqnarray*}
By (\ref{thedisk}), the first term on the right side approaches $G(\zeta_1+\bar
\zeta_2)$.  Thus it remains to show that the other terms tend to zero. Suppose
that $|\zeta_1|\le A<+\infty$.  Then,
$$\frac 1N\sum_{k=0}^N(1+{\textstyle \frac
{\zeta_1}N})^k\,\overline{q_k(1+{\textstyle \frac {\zeta_2}N})} \le 
\frac {C_1}N\sum_{k=0}^N e^{(|\zeta_1|/N -\la) k} \le \frac{C_2}N 
+ \frac {C_1}N\sum_{k=k_0}^N e^{ -\half\la k} = O\left(\frac 1N\right)\;.$$
Exactly the same estimate holds for the next term, and the last term is
clearly also $O(\frac1N)$. 

The same argument holds for the partial Bergman kernel, using Carleman's Theorem
instead of Szeg\"o's.
\end{proof}

\medskip
\subsection{The scaled zero density} We define the expected zero density function
$D^N_{\partial \Omega, \rho}$ by
$${\bf E}^N_{\partial \Omega, \rho} ( Z^N_f)=
D^N_{\partial \Omega, \rho}(z)\,\textstyle\left(\frac i2 dz\wedge d\bar
z\right)\;.$$ We recall that by Proposition \ref{circular},
$$D^N_{U,1}(z)\to \frac 1 {\pi(|z|^2-1)^2} \quad (|z|\ne 1)\;,\quad
\mbox{and}\quad\frac 1N D^N_{U,1}\to
\nu\;.$$  We now give a third limit for $D^N_{U,1}$, the scaling
limit, which we find to be universal:

\begin{prop}  Let $\Omega,\rho,\Phi$ be as above and let $\wh
D^N:=(D^N\circ\Phi\inv)|(\Phi\inv)'|^2$ be the expected zero density for the inner
product
$\langle ,\rangle _{\d\Omega,\rho}$ or for $\langle ,\rangle _\Omega$
with respect to the coordinate $\zeta=\Phi(z)$. Then
$$\frac {1}{N^2}\wh D^N\left(1+\frac \tau N\right) \to D^\infty(\tau) \qquad
\mbox{as }\ N\to \infty\;,$$
where $$ D^\infty(\tau)=\frac 1\pi (\log G)''(2\tau)={\frac {{e^{4\,\tau}}-\left
(2+4\,{\tau}^{2}\right ){e^{2\,\tau}} +1}{4\pi\left ({e^{2\,\tau}}-1\right
)^{2}{\tau}^{2}}} =\frac{1}{\pi\tau^2}-{\frac {{e^{2\,\tau}}}{4\pi\left
({e^{2\,\tau}}-1\right )^{2}}}
\;.$$  
In particular, the scaling limit density $D^\infty(\tau)$ has $1/\tau^2$ decay as
$\tau
\to
\pm \infty$.
\label{limdens}\end{prop}

\begin{proof} By Proposition \ref{poincare},
$$\wh D^N(r)=\frac 1\pi \frac{d^2}{dr^2}\log  \wh S_N(r,r)\;.$$ The conclusion
follows by substituting $r=1 +\frac \tau N$ and applying
Proposition~\ref{limszego}.\end{proof}

The graph of the scaling limit density function is given in Figure \ref{density}
below.

\begin{figure}[ht]\centering
\epsfig{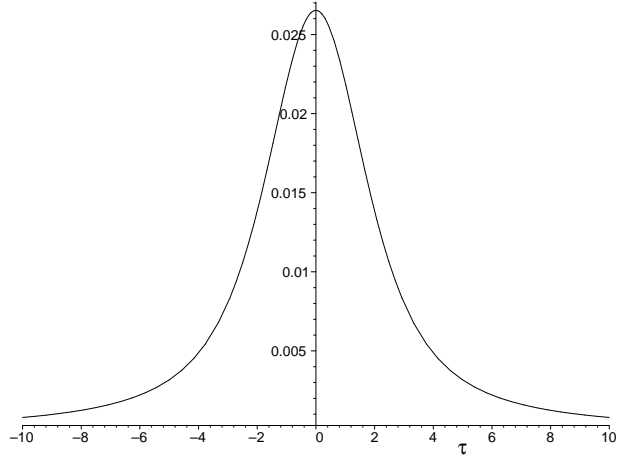}
\caption{The scaled zero density
$D^\infty(\tau)$}\label{density}
\end{figure}

\subsection{The scaling limit zero correlation functions}  Theorem~\ref{corthm} is
an immediate consequence of Propositions \ref{limszego}--\ref{limdens} and
Theorem~2.4 in \cite{BSZ1} (with $k=m=1,\ n=\ell$), which gives a universal formula for the
$\ell$-point zero correlations in terms of the projection kernel and its first and
second derivatives. (See
also \cite[Theorem~1.1]{BSZ}.) 

For the reader's convenience, we derive the formula  for the pair correlation case
$\ell=2$:  Using the coordinate $\zeta=\Phi(z)$, we define the
$2\times 2$ matrices
$${\bf A}^N(\zeta,\eta)=\left(\begin{array}{cc}\wh S_N(\zeta,\zeta)&
\wh S_N(\zeta,\eta)\\\wh S_N(\eta,\zeta) &\wh S_N(\eta,\eta)
\end{array}\right)\;,$$
$${\bf  B}^N(\zeta,\eta)= \frac{\d}{\d\bar \eta}{\bf A}^N(\zeta,\eta)\;, \qquad
{\bf C}^N(\zeta,\eta)=
\frac{\d^2}{\d \zeta\d\bar \eta}{\bf A}^N(\zeta,\eta)\;.$$ As in \cite{BSZ1,BSZ}
we let
$$\La^N={\bf  C}^N-({\bf  B}^N)^*({\bf A}^N)\inv {\bf  B}^N\;.$$ Writing $${\bf
E}^N_{\partial \Omega, \rho} ( Z^N_f\times  Z^N_f)=
\la_N(\zeta,\eta)\,d\vol_{\C^2}\;,$$ we have by  \cite[(91)]{BSZ1} or
\cite[(23)]{BSZ},
\begin{equation} \la_N = \frac{\La^N_{11} \La^N_{22}+ \La^N_{12}
\La^N_{21}}{\pi^2\det {\bf A}^N}\,.\label{lambda}\end{equation}

We consider the case of the inner product $\langle ,\rangle _{\d\Omega,\rho}$.
(The interior case is similar.) By Proposition \ref{limszego},  \begin{equation}
\label{A}\frac cN{\bf A}^N(1+\zeta_1/N, 1+\zeta_2/N)
\to
\left(\begin{array}{cc}G(\zeta_1+\bar\zeta_1)&
G(\zeta_1+\bar\zeta_2)\\G(\zeta_2+\bar\zeta_1) &G(\zeta_2+\bar\zeta_2)
\end{array}\right) \eqd {\bf A^\infty}(\zeta_1,\zeta_2)\;,\end{equation}
where $c=|\psi(z_0)|^{-2}$.
Differentiating (\ref{A}), we obtain
\begin{equation}\begin{array}{rcl}\frac {c}{N^2}{\bf B}^N(1+\zeta_1/N, 1+\zeta_2/N)
&\to&
\big(G'(\zeta_j+\bar\zeta_k)\big)_{j,k=1,2}
\eqd {\bf B^\infty}(\zeta_1,\zeta_2)
\\[10pt] \frac {c}{N^3}{\bf C}^N(1+\zeta_1/N, 1+\zeta_2/N) &\to&
\big(G''(\zeta_j+\bar\zeta_k)\big)_{j,k=1,2}
\eqd {\bf C^\infty}(\zeta_1,\zeta_2).
\end{array}\end{equation}
We write
\begin{equation}\La^\infty(\zeta_1,\zeta_2) ={\bf  C}^\infty -
({\bf  B}^\infty )^*({\bf A}^\infty )\inv {\bf  B}^\infty  =
\lim_{N\to\infty}\frac{c}{N^3}\La^N(1+\zeta_1/N, 1+\zeta_2/N)\;.
\end{equation}
Hence by (\ref{lambda}),
\begin{equation*}\frac {1}{N^4}\la_N(1+\zeta_1/N, 1+\zeta_2/N)  \to
\frac{\La^\infty _{11} \La^\infty _{22}+ \La^\infty _{12}
\La^\infty _{21}}{\pi^2\det {\bf A}^\infty }\,.\end{equation*}
Therefore by Proposition \ref{limdens},
\begin{equation}\wh K^{2 N}_{\partial \Omega, \rho}\left(1 +\frac {\zeta_1}
N,  1 +\frac {\zeta_2} N\right) \to K^{2\infty}(\zeta_1,\zeta_2)=
\frac{(\La^\infty _{11} \La^\infty _{22}+ \La^\infty _{12}
\La^\infty _{21})(\zeta_1,\zeta_2)}{D^\infty(\Re\,\zeta_1)
\,D^\infty(\Re\,\zeta_2)\,\det {\bf
A}^\infty(\zeta_1,\zeta_2) }\,.\label{formula}\end{equation}
\qed

\bigskip
The explicit expansion of (\ref{formula}) in terms of
$\zeta_1,\zeta_2$ is quite complicated. 
Using Maple$^{\rm
TM}$, we can compute the expansions of the tangential and normal
correlations
for short distances:
$$\begin{array}{l}\kappa^{\rm T}(\al)=K^{2\infty}(0,i\al)={\frac
{1}{150}}\,{\al}^{2}+{\frac {11}{42000}}\,{\al}^{4}+{
\frac {23}{5292000}}\,{\al}^{6}+{\frac {107}{5821200000}}\,{\al}
^{8} -{\frac {6659}{9535125600000}}\,{\al}^{10}\ \cdots,\\[10pt]
\kappa^\bot(\tau)=K^{2\infty}(0,\tau)={\frac {1}{150}}\,{\tau}^{2}+{\frac
{1}{1200}}\,{\tau}^{4}-{\frac
{101}{ 5292000}}\,{\tau}^{6}-{\frac {10289}{5821200000}}\,{\tau}^{8}+{\frac
{
7249481}{47675628000000}}\,{\tau}^{10}\ \cdots.\end{array}$$
Thus zeros on $\d\Omega$ repel as in the case of $\SU(2)$-polynomials
considered in \cite{Han,BSZ1}.  The behavior over long distances is
illustrated in Figures~\ref{cor-graph}--\ref{cor-tran} in the introduction. Note
that $\wh K^{\ell N}$ is invariant (only) under the $S^1$ action, and thus
$K^{2\infty}(i\al_1,i\al_2)=K^{2\infty}(0,i(\al_2-\al_1))$.  (However,
$K^{2\infty}(\tau_1,\tau_2)\neq K^{2\infty}(0,\tau_2-\tau_1)$.)

\begin{rem} The function
$$\det {\bf A}^\infty(0,i\al) = 1-\left(\frac{\sin
\,\al\!/2}{\al\!/2}\right)^2,$$
 which appears in the denominator of our
formula  for the tangential limit correlation, happens to be the limit
pair correlation function for the eigenvalues of random Hermitian matrices (see
\cite[(5.74)]{D}).  However, we are unaware of a good interpretation for $\det
{\bf A}^\infty$ in our setting.\end{rem}


\begin{thebibliography}{DPSZ}

\bibitem[Bel]{B} S. Bell, {\it The Cauchy Transform, Potential Theory, and
Conformal
Mapping}, Studies in Adv.\ Math., CRC Press, Boca Raton, FL, 1992.

\bibitem[Ber]{Berg} S. Bergman, {\it The kernel function and conformal
mapping,\/}
second revised edition, Mathematical Surveys 5, American Mathematical
Society,
Providence, RI, 1970.

\bibitem[BD]{BD} P. Bleher and X. Di,  Correlations between zeros of a
random
polynomial, {\it J. Statist.\ Phys.} 88 (1997), 269--305.

\bibitem[BR]{BR} P. Bleher and R.   Ridzal, ${\rm SU}(1,1)$ random
polynomials, {\it
J. Statist.\ Phys.} 106 (2002),  147--171.

\bibitem[BSZ1]{BSZ1} P. Bleher, B. Shiffman and S. Zelditch, Universality
and scaling of correlations between zeros on complex manifolds,  {\it
Invent.\
Math.} 142 (2000), 351--395.

\bibitem[BSZ2]{BSZ} P. Bleher, B. Shiffman, and S. Zelditch,  Correlations
between zeros and supersymmetry, {\it Comm.\ Math.\ Phys.} 224 (2001),
255--269.

\bibitem[Ca]{Ca} T. Carleman, \"Uber die Approximation analytischer
Funktionen durch
lineare Aggregate von vorgegebenen Potenzen, {\it Ark.\ Mat.\ Astr.\ Fys.}
17
(1922--23).

\bibitem[De]{D} P.  Deift, {\it Orthogonal polynomials and random matrices: a
Riemann-Hilbert approach.}
 Courant Lecture Notes in Mathematics, 3. New York University, Courant
Institute of Mathematical Sciences, New York, 1999.


\bibitem[DPSZ]{DPSZ}  A. Dembo, B. Poonen, Q.-M. Shao and  O. Zeitouni,
Random polynomials having few or no real zeros, e-print archive,
math.PR/0006113.



\bibitem[EK]{EK} A.  Edelman and E.   Kostlan, How many zeros of a random
polynomial
are real? {\it Bull.\ Amer.\ Math.\ Soc.} 32 (1995), 1--37.

\bibitem[EO]{EO} P. Erd\"os and A. C.  Offord,  On the number of real roots
of a
random algebraic equation, {\it  Proc.\ London Math.\ Soc.}  6 (1956),
139--160.

\bibitem[ET]{ET} P. Erd\"os and  P. Tur\'an,  On the distribution of roots
of
polynomials, {\it Ann.\ of Math.}  51 (1950), 105--119.

\bibitem[Ham]{H} J. M. Hammersley, The zeros of a random polynomial, {\it
Proceedings
of the Third Berkeley Symposium on Mathematical Statistics and Probability},
1954--1955, vol.\ II, pp. 89--111, University of California Press, Berkeley
and Los
Angeles, 1956.

\bibitem[Han]{Han} J. H. Hannay, Chaotic analytic zero points: exact
statistics for those
of a random spin state, {\it J. Phys.\ A\ } 29 (1996), L101--L105.

\bibitem[IM]{IM} I. A. Ibragimov and N. B.  Maslova,
The mean number of real zeros of random polynomials. I.
Coefficients with zero mean, (Russian, English summary) {\it Teor.\
Verojatnost.\ i Primenen} 16 (1971), 229--248; English translation, {\it
Theor.\ Probability
Appl.} 16 (1971), 228--248.


\bibitem[K1]{K1} M. Kac,  On the average number of real roots of a random
algebraic
equation, {\it Bull.\ Amer.\ Math.\ Soc.} 49 (1943), 314--320.

\bibitem[K2]{K2} M. Kac, On the average number of real roots of a random
algebraic
equation, II, {\it Proc.\ London Math. Soc.} 50 (1949), 390--408.

\bibitem[LO]{LO} J. E. Littlewood and A. C. Offord, On the number of real
roots of a
random algebraic equation, I, {\it  J. London Math.\ Soc.} 13  (1938),
288--295; II,
{\it  Proc.\ Cambridge Philos.\ Soc.} 35 (1939), 133--148; III, {\it Rec.\
Math.\
[Mat.\ Sbornik]} 12(54), (1943), 277--286.

\bibitem[SV]{SV} L. Shepp and R. Vanderbei, The complex zeros of random
polynomials, {\it  Trans.\ Amer.\ Math.\ Soc.} 347 (1995), 4365--4384.

\bibitem[ShZ1]{SZ} B. Shiffman and S. Zelditch,  Distribution of zeros of
random and
quantum chaotic sections of positive line bundles, {\it Comm.\ Math.\ Phys.}
200
(1999),  661--683.

\bibitem[ShZ2]{SZ2} B. Shiffman and S. Zelditch, Random polynomials with
prescribed
Newton polytope, I, e-print archive, math.AG/0203074.

\bibitem[SL]{SL} V. I. Smirnov and N. A. Lebedev,  {\it Functions of a
complex
variable: Constructive theory,\/}
M.I.T. Press, Cambridge, MA, 1968.

\bibitem[So]{So} M. Sodin,  Zeros of Gaussian analytic functions, {\it
Math.\ Res.\
Lett.} 7 (2000), 371--381.



\bibitem[Su]{Su} P. K. Suetin,  {\it Polynomials orthogonal over a region
and Bieberbach polynomials.}  Translated from the Russian by R. P. Boas.
Proceedings of the Steklov Institute of Mathematics, No. 100 (1971).
American Mathematical Society, Providence, R.I., 1974.

\bibitem[Sz1]{Sz1} G. Szeg\"o, \"Uber orthogonale Polynome, die zu einer
gegebenen
Kurve der komplexen Ebene gehoren, {\it Math.\ Zeit.} 9 (1921), 218--270.

\bibitem[Sz2]{Sz2} G.  Szeg\"o, {\it Orthogonal polynomials}, fourth
edition,
American Mathematical Society, Colloquium Publications, Vol.\ 23, American
Mathematical Society, Providence, RI, 1975.

\end{thebibliography}
\end{document}